\documentclass{amsart}
\usepackage{amsmath,amscd,amsthm,amsfonts}

\theoremstyle{plain}
\newtheorem{theorem}                 {Theorem}      [section]

\newtheorem{lemma}        [theorem]  {Lemma}

\theoremstyle{definition}
\newtheorem{example}      [theorem]  {Example}

\newtheorem{definition}   [theorem]  {Definition}

\numberwithin{equation}{section}

\begin{document}
\baselineskip 18pt \larger

\def \theo-intro#1#2 {\vskip .25cm\noindent{\bf Theorem #1\ }{\it #2}}

\newcommand{\trace}{\operatorname{trace}}

\def \dn{\mathbb D}
\def \nn{\mathbb N}
\def \zn{\mathbb Z}
\def \qn{\mathbb Q}
\def \rn{\mathbb R}
\def \cn{\mathbb C}
\def \hn{\mathbb H}
\def \P{\mathbb P}
\def \can{Ca}

\def \S{\mathcal S}
\def \A{\mathcal A}
\def \B{\mathcal B}
\def \C{\mathcal C}
\def \E{\mathcal E}
\def \F{\mathcal F}
\def \G{\mathcal G}
\def \H{\mathcal H}
\def \I{\mathcal I}
\def \L{\mathcal L}
\def \M{\mathcal M}
\def \N{\mathcal N}
\def \O{\mathcal O}
\def \P{\mathcal P}
\def \Q{\mathcal Q}
\def \R{\mathcal R}
\def \V{\mathcal V}
\def \W{\mathcal W}

\def\Re{\mathfrak R\mathfrak e}
\def\Im{\mathfrak I\mathfrak m}
\def\Co{\mathfrak C\mathfrak o}
\def\Or{\mathfrak O\mathfrak r}

\def \ip #1#2{\langle #1,#2 \rangle}
\def \spl#1#2{( #1,#2 )}

\def \lb#1#2{[#1,#2]}

\def \pror#1{\rn P^{#1}}
\def \proc#1{\cn P^{#1}}
\def \proh#1{\hn P^{#1}}

\def \gras#1#2{G_{#1}(\cn^{#2})}

\def \b{\mathfrak{b}}
\def \g{\mathfrak{g}}
\def \h{\mathfrak{h}}
\def \k{\mathfrak{k}}
\def \m{\mathfrak{m}}
\def \p{\mathfrak{p}}
\def \q{\mathfrak{q}}
\def \r{\mathfrak{r}}
\def \un{\mathfrak{u}}

\def \GLR#1{\text{\bf GL}_{#1}(\rn)}
\def \GLRP#1{\text{\bf GL}^+_{#1}(\rn)}
\def \glr#1{\mathfrak{gl}_{#1}(\rn)}
\def \GLC#1{\text{\bf GL}_{#1}(\cn)}
\def \glc#1{\mathfrak{gl}_{#1}(\cn)}
\def \GLH#1{\text{\bf GL}_{#1}(\hn)}
\def \glh#1{\mathfrak{gl}_{#1}(\hn)}
\def \GLD#1{\text{\bf GL}_{#1}(\dn)}
\def \gld#1{\mathfrak{gl}_{#1}(\dn)}

\def \SLR#1{\text{\bf SL}_{#1}(\rn)}
\def \slr#1{\mathfrak{sl}_{#1}(\rn)}
\def \SLC#1{\text{\bf SL}_{#1}(\cn)}
\def \slc#1{\mathfrak{sl}_{#1}(\cn)}

\def \O#1{\text{\bf O}(#1)}
\def \SO#1{\text{\bf SO}(#1)}
\def \so#1{\mathfrak{so}(#1)}
\def \SOs#1{\text{\bf SO}^*(#1)}
\def \sos#1{\mathfrak{so}^*(#1)}
\def \SOO#1#2{\text{\bf SO}(#1,#2)}
\def \soo#1#2{\mathfrak{so}(#1,#2)}
\def \SOC#1{\text{\bf SO}(#1,\cn)}
\def \SOc#1{\text{\bf SO}(#1,\cn)}
\def \soc#1{\mathfrak{so}(#1,\cn)}

\def \SUs#1{\text{\bf SU}^*(#1)}
\def \sus#1{\mathfrak{su}^*(#1)}

\def \U#1{\text{\bf U}(#1)}
\def \u#1{\mathfrak{u}(#1)}
\def \US#1{\text{\bf U}^*(#1)}
\def \us#1{\mathfrak{u}^*(#1)}
\def \UU#1#2{\text{\bf U}(#1,#2)}
\def \uu#1#2{\mathfrak{u}(#1,#2)}
\def \SU#1{\text{\bf SU}(#1)}
\def \su#1{\mathfrak{su}(#1)}
\def \SUU#1#2{\text{\bf SU}(#1,#2)}
\def \suu#1#2{\mathfrak{su}(#1,#2)}

\def \Sp#1{\text{\bf Sp}(#1)}
\def \sp#1{\mathfrak{sp}(#1)}
\def \Spp#1#2{\text{\bf Sp}(#1,#2)}
\def \spp#1#2{\mathfrak{sp}(#1,#2)}
\def \SpR#1{\text{\bf Sp}(#1,\rn)}
\def \spR#1{\mathfrak{sp}(#1,\rn)}
\def \SpC#1{\text{\bf Sp}(#1,\cn)}
\def \spc#1{\mathfrak{sp}(#1,\cn)}

\def \d#1{\mathfrak{d}(#1)}
\def \s#1{\mathfrak{s}(#1)}
\def \sym#1{\text{Sym}(\rn^{#1})}
\def \symc#1{\text{Sym}(\cn^{#1})}

\def \gradh#1{\text{grad}_{\H}(#1 )}
\def \gradv#1{\text{grad}_{\V}(#1 )}

\def \nab#1#2{\hbox{$\nabla$\kern -.3em\lower 1.0 ex
    \hbox{$#1$}\kern -.1 em {$#2$}}}

\allowdisplaybreaks

\title{Harmonic morphisms from the classical
\\compact semisimple Lie groups\\(version 1.059)}

\author{Sigmundur Gudmundsson}
\author{Anna Sakovich}



\address
{Mathematics, Faculty of Science, Lund University, Box 118, S-221
00 Lund, Sweden} \email{Sigmundur.Gudmundsson@math.lu.se}

\address
{Faculty of Pre-University Preparation, Belarusian State
University, Oktyabrskaya Str. 4, Minsk 220030, Belarus}
\email{Anya\_Sakovich@tut.by}

\begin{abstract}
In this paper we introduce a new method for manufacturing harmonic
morphisms from semi-Riemannian manifolds.  This is employed to
yield a variety of new examples from the compact Lie groups $\SO
n$, $\SU n$ and $\Sp n$ equipped with their standard Riemannian
metrics. We develop a duality principle and show how this can be
used to construct the first known examples of harmonic morphisms
from the non-compact Lie groups $\SLR n$, $\SUs{2n}$, $\SpR n$,
$\SOs{2n}$, $\SOO pq$, $\SUU pq$ and $\Spp pq$ equipped with their
standard dual semi-Riemannian metrics.
\end{abstract}

\maketitle

\section{Introduction}

The notion of a minimal submanifold of a given ambient space is of
great importance in differential geometry. Harmonic morphisms
$\phi:(M,g)\to(N,h)$ between semi-Riemannian manifolds are useful
tools for the construction of such objects. They are solutions to
over-determined non-linear systems of partial differential
equations determined by the geometric data of the manifolds
involved. For this reason harmonic morphisms are difficult to find
and have no general existence theory, not even locally.

For the existence of harmonic morphisms $\phi:(M,g)\to (N,h)$ it
is an advantage that the target manifold $N$ is a surface i.e. of
dimension $2$. In this case the problem is invariant under
conformal changes of the metric on $N^2$. Therefore, at least for
local studies, the codomain can be taken to be the complex plane
with its standard flat metric.   For the general theory of
harmonic morphisms between semi-Riemannian manifolds we refer to
the excellent book \cite{Bai-Woo-book} and the regularly updated
on-line bibliography \cite{Gud-bib}.

The Riemannian manifold Sol is one of the eight famous
3-dimensional homogeneous geometries and has the structure of a
Lie group compatible with its Riemannian metric.  Baird and Wood
have shown in \cite{Bai-Woo-1} that Sol does not allow any local
harmonic morphisms with values in a surface. This fact has been
the main motivation for the research leading to this paper.

We introduce the notion of an {\it eigenfamily} of complex valued
functions on a given semi-Riemannian manifold $(M,g)$. We show how
such a family can be used to construct a variety of harmonic
morphisms on open and dense subsets of $M$.

We focus our attention on the classical semi-simple Lie groups and
construct eigenfamilies in the compact Riemannian cases of
$$\SO n,\ \ \SU n\ \ \text{and}\ \ \Sp n$$ inducing a variety of
new harmonic morphisms on these important spaces.  The examples
constructed by our new method are the first harmonic morphisms on
$\SO n$, $\SU n$ and $\Sp n$ which are not invariant under the
action of the subgroups
$$\SO p\times\SO q,\ \ \text{\bf S}(\U p\times\U q)\ \ \text{and}
\ \ \Sp p\times\Sp q,$$ respectively, with $n=p+q$.  For the known
invariant solutions leading to harmonic morphisms on the
Grassmannians, see \cite{Gud-Sve-1}.

In Theorem \ref{theo:duality} we prove a useful duality principle
and show how this can be used to obtain eigenfamilies on the
non-compact semi-Riemannian Lie groups
$$\SLR n,\ \ \SUs{2n},\ \ \SpR n,$$
$$\SOs{2n},\ \ \SOO pq,\ \ \SUU pq\ \ \text{and}\ \ \Spp pq.$$
This leads to the construction of the first known examples of
harmonic morphisms in all these cases.   It should be noted that
the non-compact semi-simple Lie groups
$$\SOC n,\ \ \SLC n\ \ \text{and}\ \ \SpC n$$ are complex manifolds and
hence their coordinate functions form orthogonal harmonic
families, see Definition \ref{defi:eigen}.  This means that in
these cases the problem is more or less trivial.

Throughout this article we assume, when not stating otherwise,
that all our objects such as manifolds,  maps etc. are smooth,
i.e. in the $C^{\infty}$-category. For our notation concerning Lie
groups we refer to the wonderful book \cite{Kna}.

\section{Harmonic Morphisms}

Let $M$ and $N$ be two manifolds of dimensions $m$ and $n$,
respectively. Then a semi-Riemannian metric $g$ on $M$ gives rise
to the notion of a Laplacian on $(M,g)$ and real-valued harmonic
functions $f:(M,g)\to\rn$. This can be generalized to the concept
of a harmonic map $\phi:(M,g)\to (N,h)$ between semi-Riemannian
manifolds being a solution to a semi-linear system of partial
differential equations, see \cite{Bai-Woo-book}.

\begin{definition}
A map $\phi:(M,g)\to (N,h)$ between semi-Riemannian manifolds is
called a {\it harmonic morphism} if, for any harmonic function
$f:U\to\rn$ defined on an open subset $U$ of $N$ with
$\phi^{-1}(U)$ non-empty, the composition
$f\circ\phi:\phi^{-1}(U)\to\rn$ is a harmonic function.
\end{definition}

The following characterization of harmonic morphisms between
semi-Riemannian manifolds is due to Fuglede and generalizes the
corresponding well-known result of \cite{Fug-1,Ish} in the
Riemannian case.  For the definition of horizontal conformality we
refer to \cite{Bai-Woo-book}.

\begin{theorem}\cite{Fug-2}
  A map $\phi:(M,g)\to (N,h)$ between semi-Rie\-mannian manifolds is a
  harmonic morphism if and only if it is a horizontally (weakly)
  conformal harmonic map.
\end{theorem}

The next result generalizes the corresponding well-known theorem
of Baird and Eells in the Riemannian case, see \cite{Bai-Eel}. It
gives the theory of harmonic morphisms a strong geometric flavour
and shows that the case when $n=2$ is particularly interesting. In
that case the conditions characterizing harmonic morphisms are
independent of conformal changes of the metric on the surface
$N^2$.  For the definition of horizontal homothety we refer to
\cite{Bai-Woo-book}.

\begin{theorem}\cite{Gud-1}\label{theo:semi-B-E}
Let $\phi:(M,g)\to (N^n,h)$ be a horizontally conformal submersion
from a semi-Riemannian manifold $(M,g)$ to a Riemannian manifold
$(N,h)$. If
\begin{enumerate}
\item[i.] $n=2$ then $\phi$ is harmonic if and only if $\phi$ has
minimal fibres, \item[ii.] $n\ge 3$ then two of the following
conditions imply the other:
\begin{enumerate}
\item $\phi$ is a harmonic map, \item $\phi$ has minimal fibres,
\item $\phi$ is horizontally homothetic.
\end{enumerate}
\end{enumerate}
\end{theorem}

In what follows we are mainly interested in complex valued
functions $\phi,\psi:(M,g)\to\cn$ from semi-Riemannian manifolds.
In this situation the metric $g$ induces the complex-valued
Laplacian $\tau(\phi)$ and the gradient $\text{grad}(\phi)$ with
values in the complexified tangent bundle $T^{\cn}M$ of $M$.  We
extend the metric $g$ to be complex bilinear on $T^{\cn} M$ and
define the symmetric bilinear operator $\kappa$ by
$$\kappa(\phi,\psi)= g(\text{grad}(\phi),\text{grad}(\psi)).$$ Two
maps $\phi,\psi: M\to\cn$ are said to be {\it orthogonal} if
$$\kappa(\phi,\psi)=0.$$  The harmonicity and horizontal
conformality of $\phi:(M,g)\to\cn$ are given by the following
relations
$$\tau(\phi)=0\ \ \text{and}\ \ \kappa(\phi,\phi)=0.$$

\begin{definition}\label{defi:eigen}
Let $(M,g)$ be a semi-Riemannian manifold.  Then a set
$$\E =\{\phi_i:M\to\cn\ |\ i\in I\}$$ of complex valued functions
is said to be an {\it eigenfamily} on $M$ if there exist complex
numbers $\lambda,\mu\in\cn$ such that
$$\tau(\phi)=\lambda\phi\ \ \text{and}\ \ \kappa
(\phi,\psi)=\mu\phi\psi$$ for all $\phi,\psi\in\E$. A set
$$\Omega =\{\phi_i:M\to\cn\ |\ i\in I\}$$ is said to be an {\it
orthogonal harmonic} family on $M$ if  for all
$\phi,\psi\in\Omega$
$$\tau(\phi)=0\ \ \text{and}\ \ \kappa(\phi,\psi)=0.$$
\end{definition}

The next result shows that an eigenfamily on a semi-Riemannian
manifold can be used to produce a variety of local harmonic
morphisms.

\begin{theorem}\label{theo:rational}
Let $(M,g)$ be a semi-Riemannian manifold and $$\E =\{\phi_1,\dots
,\phi_n\}$$ be a finite eigenfamily of complex valued functions on
$M$. If $P,Q:\cn^n\to\cn$ are linearily independent homogeneous
polynomials of the same positive degree then the quotient
$$\frac{P(\phi_1,\dots ,\phi_n)}{Q(\phi_1,\dots ,\phi_n)}$$ is a
non-constant harmonic morphism on the open and dense subset
$$\{p\in M| \ Q(\phi_1(p),\dots ,\phi_n(p))\neq 0\}.$$
\end{theorem}

A proof of Theorem \ref{theo:rational} can be found in Appendix
\ref{app:general}. For orthogonal harmonic families we have the
following useful result.

\begin{theorem}\cite{Gud-1}\label{theo:orthogonal}
Let $(M,g)$ be a semi-Riemannian manifold and
$$\{\phi_k:M\to\cn\ |\ k=1,\dots ,n\}$$ be a finite orthogonal
harmonic family on $(M,g)$.  Let $\Phi:M\to\cn^n$ be the map given
by $\Phi=(\phi_1,\dots,\phi_n)$ and $U$ be an open subset of
$\cn^n$ containing the image $\Phi(M)$ of $\Phi$. If
$$\H=\{h_i:U\to\cn\ |\ i\in I\}$$ is a family of holomorphic
functions then $$\{\psi:M\to\cn\ |\ \psi=h(\phi_1,\dots ,\phi_n
),\ h\in\H\}$$ is an orthogonal harmonic family on $M$.
\end{theorem}

\section{The Riemannian Lie group $\GLC n$}

Let $G$ be a Lie group with Lie algebra $\g$ of left-invariant
vector fields on $G$.  Then a Euclidean scalar product $g$ on the
algebra $\g$ induces a left-invariant Riemannian metric on the
group $G$ and turns it into a Riemannian manifold. If $Z$ is a
left-invariant vector field on $G$ and $\phi,\psi:U\to\cn$ are two
complex valued functions defined locally on $G$ then the first and
second order derivatives satisfy
$$Z(\phi)(p)=\frac {d}{ds}[\phi(p\cdot\exp(sZ))]\big|_{s=0},$$
$$Z^2(\phi)(p)=\frac {d^2}{ds^2}[\phi(p\cdot\exp(sZ))]\big|_{s=0}.$$
The tension field $\tau(\phi)$ and the $\kappa$-operator
$\kappa(\phi,\psi)$ are given by
$$\tau(\phi)=\sum_{Z\in\B}Z^2(\phi)\ \ \text{and}\ \
\kappa(\phi,\psi)=\sum_{Z\in\B}Z(\phi)Z(\psi)$$ where $\B$ is any
orthonormal basis of the Lie algebra $\g$.

Let $\GLC n$ be the complex general linear group equipped with its
standard Riemannian metric induced by the Euclidean scalar product
on the Lie algebra $\glc n$ given by
$$g(Z,W)=\Re\trace ZW^*.$$  For $1\le i,j\le n$ we shall by
$E_{ij}$ denote the element of $\glr n$ satisfying
$$(E_{ij})_{kl}=\delta_{ik}\delta_{jl}$$ and by $D_t$ the diagonal
matrices $$D_t=E_{tt}.$$ For $1\le r<s\le n$ let $X_{rs}$ and
$Y_{rs}$ be the matrices satisfying
$$X_{rs}=\frac 1{\sqrt 2}(E_{rs}+E_{sr}),\ \ Y_{rs}=\frac
1{\sqrt 2}(E_{rs}-E_{sr}).$$ With the above notation we have the
following easily verified matrix identities
$$\sum_{r<s}X_{rs}^2=\frac {(n-1)}2I_n,\ \ \
\sum_{r<s}Y_{rs}^2=-\frac {(n-1)}2I_n,\ \ \
\sum_{t=1}^nD_t^2=I_n,$$

$$\sum_{r<s}X_{rs}E_{jl}X^t_{rs}=\frac 12(E_{lj}+\delta
_{lj}(I_n-2E_{lj})),$$

$$\sum_{r<s}Y_{rs}E_{jl}Y^t_{rs}=-\frac 12( E _{lj}-\delta _{lj}I_n),$$

$$
\sum_{t=1}^nD_{t}E_{jl}D^t_{t}=\  \delta _{jl} E _{lj}.
$$

\section{The Riemannian Lie group $\SO n$}

In this section we construct eigenfamilies of complex valued
functions on the special orthogonal group
$$\SO{n}=\{x\in\GLR{n}\ |\ x\cdot x^t=I_n,\ \det x =1\}.$$ The Lie
algebra $\so n$ of $\SO n$ is the set of skew-symmetric matrices
$$\so n=\{X\in\glr n |\ X+X^t=0\}$$ and for this we have the
canonical orthonormal basis
$$\{Y_{rs}|\ 1\le r<s\le n\}.$$

\begin{lemma}\label{lemm:real}
For $1\le i,j\le n$ let $x_{ij}:\SO n\to\rn$ be the real valued
coordinate functions given by
$$x_{ij}:x\mapsto e_i\cdot x\cdot e_j^t$$ where $\{e_1,\dots
,e_n\}$ is the canonical basis for $\rn^n$. Then the following
relations hold  $$\tau(x_{ij})=-\frac {(n-1)}2x_{ij},$$
$$\kappa(x_{ij},x_{kl})=-\frac 12(x_{il}x_{kj}-\delta_{jl}
\sum_{t=1}^nx_{it}x_{kt}).$$
\end{lemma}

\begin{proof}  It follows directly from the definition of the functions
$x_{ij}$ that if $X$ is an element of the Lie algebra $\so n$ then
the first and second order derivatives satisfy
$$X(x_{ij}):x\mapsto e_i\cdot x\cdot X\cdot e_j^t\ \ \text{and}
\ \  X^2(x_{ij}):x\mapsto e_i\cdot x\cdot X^2\cdot e_j^t.$$
Employing the above mentioned matrix identities we then yield
$$\tau(x_{ij})=\sum_{r<s}Y_{rs}^2(x_{ij})=
\sum_{r<s}e_i\cdot x\cdot Y_{rs}^2\cdot e_j^t=-\frac
{(n-1)}2x_{ij},$$
\begin{eqnarray*}
\kappa(x_{ij},x_{kl}) &=&\sum_{r<s}e_i\cdot x\cdot Y_{rs}\cdot e_j^t
\cdot e_l\cdot Y_{rs}^t\cdot x^t\cdot e_k^t\\
&=&e_i\cdot x\cdot\big(\sum_{r<s}Y_{rs}\cdot E_{jl}\cdot
Y_{rs}^t\big)\cdot x^t\cdot e_k^t\\
&=&-\frac 12(x_{il}x_{kj}-\delta_{jl}\sum_{t=1}^nx_{it}x_{kt}).
\end{eqnarray*}
\end{proof}

Let $P,Q:\SO n\to\cn$ be homogeneous polynomials of the coordinate
functions $x_{ij}:\SO n\to\cn$ of degree one i.e. of the form
$$P=\trace (A\cdot x^t)=\sum_{i,j=1}^n a_{ij}x_{ij}\ \ \text{and}
\ \ Q=\trace (B\cdot x^t)=\sum_{k,l=1}^n b_{kl}x_{kl}$$ for some
$A,B\in\cn^{n\times n}$.  As a direct consequence of Lemma
\ref{lemm:real} we then yield
\begin{eqnarray*}
& &PQ+2\kappa(P,Q)\\
&=&\sum_{i,j,k,l=1}^na_{ij}b_{kl}x_{ij}x_{kl}
+2\sum_{i,j,k,l=1}^na_{ij}b_{kl}\kappa(x_{ij},x_{kl})\\
&=&\sum_{i,j,k,l=1}^na_{ij}b_{kl}x_{ij}x_{kl}
-\sum_{i,j,k,l=1}^na_{ij}b_{kl}x_{kj}x_{il}
+\sum_{i,j,k,t=1}^na_{ij}b_{kj}x_{it}x_{kt}\\
&=&\sum_{i,j,k,l=1}^n(a_{ij}b_{kl}-a_{kj}b_{il})x_{ij}x_{kl}
+\trace (AB^txx^t).
\end{eqnarray*}
Comparing coefficients we see that $PQ+2\kappa (P,Q)=0$ if
$AB^t=0$ and
$$\det\begin{pmatrix} a_{ij}& b_{il} \\a_{kj} &
b_{kl}\end{pmatrix} =(a_{ij}b_{kl}-a_{kj}b_{il})=0$$ for all $1\le
i,j,k,l\le n$.

\begin{theorem}\label{theo:real}
Let $V$ be a maximal isotropic subspace of $\cn^n$ and $p\in\cn^n$
be a non-zero element.  Then the set
$$\E_V(p)=\{\phi_a:\SO n\to\cn\ |\ \phi_a(x)=\trace (p^tax^t),\
a\in V\}$$ of complex valued functions is an eigenfamily on $\SO
n$.
\end{theorem}

\begin{proof}  Assume that $a,b\in V$ and define $A=p^ta$ and
$B=p^tb$.  By construction any two columns of the matrices $A$ and
$B$ are linearly dependent.  This means that for all $1\le
i,j,k,l\le n$
$$\det\begin{pmatrix} a_{ij}& b_{il}
\\a_{kj} & b_{kl}\end{pmatrix} =(a_{ij}b_{kl}-a_{kj}b_{il})=0.$$
Furthermore we have $AB^t=0$. Hence $P^2+2\kappa (P,P)=0$,
$PQ+2\kappa (P,Q)=0$, $Q^2+2\kappa (Q,Q)=0$ and the statement
follows directly from Lemma \ref{lemm:real} and the calculations
above.
\end{proof}

Applying the fact that $xx^t=I$ for each $x\in\SO n$ we get a
simplified formula for the $\kappa$ operator
$$\kappa(x_{ij},x_{kl})=\frac
12(\delta_{ik}\delta_{jl}-x_{il}x_{kj}).$$  With a similar
analysis to that above one yields the result of Theorem
\ref{theo:real-special}.  It should be noted that having employed
the special property $xx^t=I$ the eigenfamily $\E(p)$ can not be
used directly in the duality of Theorem \ref{theo:duality}.

\begin{theorem}\label{theo:real-special}
Let $p$ be a non-zero isotropic element of $\cn^n$ i.e. such that
$(p,p)=0$. Then the set
$$\E(p)=\{\phi_a:\SO n\to\cn\ |\ \phi_a(x)=\trace (p^tax^t),\
a\in \cn^n\}$$ of complex valued functions is an eigenfamily on
$\SO n$.
\end{theorem}

\begin{example}
For $z,w\in\cn$ let $p$ be the isotropic element of the
$4$-dimensional complex vector space $\cn^4$ given by
$$p(z,w)=(1+zw,i(1-zw),i(z+w),z-w).$$   This gives us the
complex $2$-dimensional deformation of eigenfamilies $\E_p$ each
consisting of complex valued functions
$$\phi_a:\SO 4\to\cn$$ of the form
\begin{eqnarray*}
\phi_a(x)=(1+zw)(a_1x_{11}+a_2x_{21}+a_3x_{31}+a_4x_{41})& &\\
         +i(1-zw)(a_1x_{12}+a_2x_{22}+a_3x_{32}+a_4x_{42})& &\\
         +i(z+w)(a_1x_{13}+a_2x_{23}+a_3x_{33}+a_4x_{43})& &\\
         +(z-w)(a_1x_{14}+a_2x_{24}+a_3x_{34}+a_4x_{44}).& &
\end{eqnarray*}
\end{example}

\section{The Riemannian Lie group $\SU n$}

In this section we construct eigenfamilies of complex valued
functions on the unitary group $\U n$.  They can be used to
construct local harmonic morphisms on the special unitary group
$\SU n$. The unitary group $\U n$ is the compact subgroup of $\GLC
n$ given by $$\U{n}=\{z\in\GLC{n}|\ z\cdot z^*=I_n\}.$$  The
circle group $S^1=\{w\in\cn | \ |w|=1\}$ acts on the unitary group
$\U n$ by multiplication $(e^{i\theta},z)\mapsto e^{i\theta} z$
and the orbit space of this action is the special unitary group
$$\SU n=\{ z\in\U n|\ \det z = 1\}.$$ The Lie algebra $\u n$ of the
unitary group $\U n$ satisfies
$$\u{n}=\{Z\in\cn^{n\times n}|\ Z+Z^*=0\}$$ and for this
we have the canonical orthonormal basis
$$\{Y_{rs}, iX_{rs}|\ 1\le r<s\le n\}\cup\{iD_t|\ t=1,\dots ,n\}.$$

\begin{lemma}\label{lemm:complex}
For $1\le i,j\le n$ let $z_{ij}:\U n\to\cn$ be the complex valued
coordinate functions given by
$$z_{ij}:z\mapsto e_i\cdot z\cdot e_j^t$$ where $\{e_1,\dots
,e_n\}$ is the canonical basis for $\cn^n$. Then the following
relations hold $$\tau(z_{ij})= -n z_{ij}\ \ \text{and}\ \
\kappa(z_{ij},z_{kl})= -z_{il}z_{kj}.$$
\end{lemma}

\begin{proof}  The proof is similar to that of Lemma
\ref{lemm:real}.
\end{proof}

Let $P,Q:\U n\to\cn$ be homogeneous polynomials of the coordinate
functions $z_{ij}:\U n\to\cn$ of degree one i.e. of the form
$$P=\trace (A\cdot z^t)=\sum_{i,j=1}^n a_{ij}z_{ij}\ \ \text{and}
\ \ Q=\trace (B\cdot z^t)=\sum_{k,l=1}^n b_{kl}z_{kl}$$ for some
$A,B\in\cn^{n\times n}$.  As a direct consequence of Lemma
\ref{lemm:complex} we then yield
$$PQ+\kappa(P,Q)=\sum_{i,j,k,l=1}^n(a_{ij}b_{kl}-a_{kj}b_{il})z_{ij}z_{kl}.$$
Comparing coefficients we see that $\kappa (P,Q)+PQ=0$ if for all
$1\le i,j,k,l\le n$ $$\det\begin{pmatrix} a_{ij}& b_{il}
\\a_{kj} & b_{kl}\end{pmatrix} =(a_{ij}b_{kl}-a_{kj}b_{il})=0.$$

\begin{theorem}\label{theo:complex}
Let $p$ be a non-zero element of $\cn^n$. Then the set
$$\E(p)=\{\phi_a:\U n\to\cn\ |\ \phi_a(z)=\trace (p^taz^t),\
a\in \cn^n\}$$ of complex valued functions is an eigenfamily on
$\U n$.
\end{theorem}

\begin{proof} Assume that $a,b\in\cn^n$ and define $A=p^ta$ and
$B=p^tb$. By construction any two columns of the matrices $A$ and
$B$ are linearly dependent.  This means that for all $1\le
i,j,k,l\le n$
$$\det\begin{pmatrix} a_{ij}& b_{il}
\\a_{kj} & b_{kl}\end{pmatrix} =(a_{ij}b_{kl}-a_{kj}b_{il})=0$$
so $P^2+\kappa (P,P)=0$, $PQ+\kappa (P,Q)=0$ and $Q^2+\kappa
(Q,Q)=0$. The statement now follows directly from Lemma
\ref{lemm:complex}.
\end{proof}

It should be noted that the local harmonic morphisms on the
unitary group $\U n$ that we obtain by employing Theorem
\ref{theo:rational} are invariant under the cirle action and hence
induce local harmonic morphisms on the special unitary group $\SU
n$.

\begin{example}
The $3$-dimensional sphere $S^3$ is diffeomorphic to the special
unitary group $\SU 2$ given by $$\SU 2=\{\begin{pmatrix} z& w
\\-\bar w & \bar z\end{pmatrix} |\ |z|^2+|w|^2=1\}.$$
For $p=(1,0)\in\cn^2$ we get the eigenfamily $$\E(p)=\{\phi_a:\U
n\to\cn\ |\ \phi_a(z)=a_1z+a_2w,\ a=(a_1,a_2)\in \cn^n\}.$$ By
choosing $a=(1,0)$ and $b=(0,1)$ and applying Theorem
\ref{theo:rational} we obtain the well known globally defined
harmonic morphism
$$\phi=\frac{\phi_a}{\phi_b}:\SU 2\to S^2$$ called the Hopf map satisfying
$$\phi(\begin{pmatrix} z& w \\-\bar w & \bar z\end{pmatrix})=\frac zw.$$
\end{example}

\section{The Riemannian Lie group $\Sp n$}
In this section we construct eigenfamilies of complex valued
functions from the quaternionic unitary group $\Sp n$ being the
intersection of the unitary group $\U{2n}$ and the standard
representation of the quaternionic general linear group $\GLH n$
in $\cn^{2n\times 2n}$ given by
$$(z+jw)\mapsto q=\begin{pmatrix}z & w \\ -\bar w & \bar
z\end{pmatrix}.$$ The Lie algebra $\sp n$ of $\Sp n$ satisfies
$$\sp{n}=\{\begin{pmatrix} Z & W
\\ -\bar W & \bar Z\end{pmatrix}\in\cn^{2n\times 2n}
\ |\ Z^*+Z=0,\ W^t-W=0\}$$ and for this we have the standard
orthonormal basis which is the union of the following three sets
$$\{\frac 1{\sqrt 2}\begin{pmatrix}Y_{rs} & 0 \\
0 & Y_{rs}\end{pmatrix},\frac 1{\sqrt 2}\begin{pmatrix}iX_{rs} & 0 \\
0 & -iX_{rs}\end{pmatrix}\ |\ 1\le r<s\le n\},$$

$$\{\frac 1{\sqrt 2}\begin{pmatrix}0 & iX_{rs} \\
iX_{rs} & 0\end{pmatrix},\frac 1{\sqrt 2}\begin{pmatrix}0 & X_{rs} \\
-X_{rs} & 0\end{pmatrix}\ |\ 1\le r<s\le n\},$$

$$\{\frac 1{\sqrt 2}\begin{pmatrix}iD_{t} & 0 \\
0 & -iD_{t}\end{pmatrix},\frac 1{\sqrt 2}\begin{pmatrix}0 & iD_{t}  \\
iD_{t} & 0\end{pmatrix},\frac 1{\sqrt 2}\begin{pmatrix}0 & D_{t}  \\
-D_{t} & 0\end{pmatrix}\ |\ 1\le t\le n\}.$$

\begin{lemma}\label{lemm:quaternionic}
For $1\le i,j\le n$ let $z_{ij},w_{ij}:\Sp n\to\cn$ be the complex
valued coordinate functions given by
$$z_{ij}:g\mapsto e_i\cdot g\cdot e_j^t,\ \
w_{ij}:g\mapsto e_i\cdot g\cdot e_{n+j}^t$$ where $\{e_1,\dots
,e_{2n}\}$ is the canonical basis for $\cn^{2n}$. Then the
following relations hold
$$\tau(z_{ij})= -\frac{2n+1}2\cdot z_{ij},\ \ \tau(w_{ij})=
-\frac{2n+1}2\cdot w_{ij},$$

$$\kappa(z_{ij},z_{kl})=-\frac 12\cdot z_{il}z_{kj},\ \
\kappa(w_{ij},w_{kl})=-\frac 12\cdot w_{il}w_{kj},$$

$$\kappa(z_{ij},w_{kl})=-\frac 12\big[w_{il}z_{kj}
-{\delta_{jl}}\cdot \sum_{t=1}^n(z_{it}w_{kt}
-w_{it}z_{kt})\big].$$
\end{lemma}

\begin{proof}  The proof is similar to that of Lemma
  \ref{lemm:real} but more involved.
\end{proof}

\begin{theorem}\label{theo:quaternionic-special}
Let $p$ be a non-zero element of $\cn^n$. Then the set
$$\E(p)=\{\phi_{ab}:\Sp n\to\cn\ |\ \phi_{ab}(g)=\trace (p^taz^t+p^tbw^t),\
a,b\in \cn^n\}$$ of complex valued functions is an eigenfamily on
$\Sp n$.
\end{theorem}

\begin{proof} Let $a,b,c,d$ be arbitrary elements of $\cn^n$ and define the
complex valued functions $P,Q:\Sp n\to\cn$ by
$$P=\trace (p^taz^t+p^tbw^t)\ \
\text{and}\ \ Q=\trace (p^tcz^t+p^tdw^t).$$ Then a simple
calculation shows that
$$PQ+2\kappa(P,Q)=[(a,d)-(b,c)]
\sum_{i,k,t=1}^n(z_{it}w_{kt}-w_{it}z_{kt})=0.$$
Automatically we also get $P^2+2\kappa(P,P)=0$ and $Q^2+2\kappa(Q,Q)=0$.
\end{proof}

\section{The Duality}

In this section we show how a real analytic eigenfamily $\E$ on a
semi-Riemannian non-compact semi-simple Lie group $G$ gives rise
to a real-analytic eigenfamily $\E^*$ on its Riemannian compact
dual $U$ and vice versa.  The method of proof is borrowed from a
related duality principle for harmonic morphisms from Riemannian
symmetric spaces, see \cite{Gud-Sve-1}.

Let $W$ be an open subset of $G$ and $\phi:W\to\cn$ be a real
analytic map.  Let $G^\cn$ denote the complexification of the Lie
group $G$. Then $\phi$ extends uniquely to a holomorphic map
$\phi^\cn:W^\cn\to\cn$ from some open subset $W^\cn$ of $G^\cn$.
By restricting this map to $U\cap W^\cn$ we obtain a real analytic
map $\phi^*:W^*\to\cn$ from some open subset $W^*$ of $U$.

\begin{theorem}\label{theo:duality}
Let $\E$ be a family of maps $\phi:W\to\cn$ and $\E^*$ be the dual
family consisting of the maps $\phi^*:W^*\to\cn$ constructed as
above.  Then $\E^*$ is an eigenfamily if and only if $\E$ is an
eigenfamily.
\end{theorem}

\begin{proof}  Let $\g=\k+\p$ be a Cartan decomposition of the Lie
algebra of $G$ where $\k$ is the Lie algebra of a maximal compact
subgroup $K$.  Furthermore let the left-invariant vector fields
$X_1,\dots,X_n\in\p$ form a global orthonormal frame for the
distribution generated by $\p$ and similarly $Y_1,\dots,Y_m\in\k$
form a global orthonormal frame for the distribution generated by
$\k$. We shall now assume that $\phi$ and $\psi$ are elements of
the eigenfamily $\E$ on the semi-Riemannian $W$ i.e.
$$\tau (\phi)=-\sum_{k=1}^m Y_k^{2}(\phi) +\sum_{k=1}^n
X_k^{2}(\phi)=\lambda\cdot\phi,$$ $$\kappa
(\phi,\psi)=-\sum_{k=1}^m Y_k(\phi)Y_k(\psi)+\sum_{k=1}^n
X_k(\phi)X_k(\psi)=\mu\cdot\phi\cdot\psi.$$ By construction and by
the unique continuation property of real analytic functions the
extension ${\phi}^\cn$ of $\phi$ satisfies the same equations.

The Lie algebra of $U$ has the decomposition $\un=\k+i\p$ and the
left-invariant vector fields $iX_1,\dots,iX_n\in\ i\p$ form a
global orthonormal frame for the distribution generated by $i\p$.
Then

$$\tau (\phi^*)=\sum_{k=1}^m Y_k^{2}(\phi^*)
+\sum_{k=1}^n (iX_k)^{2}(\phi^*)=-\lambda\cdot\phi^*$$

$$\kappa (\phi^*,\phi^*)=\sum_{k=1}^m Y_k(\phi^*)Y_k(\psi^*)
+\sum_{k=1}^n
(iX_k)(\phi^*)(iX_k)(\psi^*)=-\mu\cdot\phi^*\cdot\psi^*.$$ This
shows that $\E^*$ is an eigenfamily. The converse is similar.
\end{proof}

\section{The semi-Riemannian Lie group $\GLC n$}

Let $h$ be the standard left-invariant semi-Riemannian metric on
the general linear group $\GLC n$ induced by the semi-Euclidean
scalar product on the Lie algebra $\glc n$ given by
$$h(Z,W)=\Re\trace ZW.$$ Then we have the orthogonal
decomposition $$\glc n=\W_+\oplus\W_-$$ of Lie algebra $\glc n$
where
$$\W_+=\{Z\in\glc n |\ Z-Z^*=0\}$$ is the subspace of Hermitian
matrices and
$$\W_-=\{Z\in\glc n |\ Z+Z^*=0\}$$ is the subspace of skew-Hermitian
matrices.  The scalar product is positive definite on $\W_+$ and
negative definite on $\W_-$.  This means that for two complex
valued functions $\phi,\psi:U\to\cn$ locally defined on $\GLC n$
the differential operators $\tau$ and the $\kappa$ satisfy
$$\tau(\phi)=\sum_{Z\in\B_+}Z^2(\phi)-\sum_{Z\in\B_-}Z^2(\phi),$$
$$\kappa(\phi,\psi)=\sum_{Z\in\B_+}Z(\phi)Z(\psi)
-\sum_{Z\in\B_-}Z(\phi)Z(\psi)$$ where $\B_+$ and $\B_-$ are
orthonormal bases for $\W_+$ and $\W_-$, respectively.

Employing the duality principle of Theorem \ref{theo:duality} we
can now easily construct harmonic morphisms from the non-compact
semi-Riemannian Lie groups
$$\SLR n,\ \SUs{2n},\ \SpR n,\ \SOs{2n},\ \SOO pq,\ \SUU pq,\ \Spp pq$$
via the following classical dualities $G\cong U$:

$$\SLR n=\{x\in\GLR{n}\ |\ \det x=1\} \cong \SU n,$$

$$\SUs{2n}=\{g=\begin{pmatrix} z & w
\\ -\bar w & \bar z\end{pmatrix}
\ |\ g\in\SLC {2n}\}\cong \SU{2n},$$

$$\SpR n =\{g\in\SLR{2n}\ |\ g\cdot J_{n}\cdot
g^t=J_{n}\} \cong\Sp n,$$

$$\SOs{2n}=\{z\in\SUU nn\ |\ z\cdot I_{nn}\cdot J_{n}\cdot
z^t=I_{nn}\cdot J_{n}\}\cong\SO{2n},$$

$$\SOO pq=\{x\in\SLR{p+q}\ |\ x\cdot I_{pq}\cdot
x^t=I_{pq}\}\cong\SO{p+q},$$

$$\SUU pq=\{z\in\SLC{p+q}\ |\ z\cdot I_{pq}\cdot
z^*=I_{pq}\}\cong\SU{p+q},$$

$$\Spp pq=\{g\in\GLH{p+q}\ |\ g\cdot I_{pq}\cdot g^*
=I_{pq}\}\cong \Sp{p+q}.$$ Here we have used the standard notation
$$I_{pq}=\begin{pmatrix} -I_p& 0
\\0 & I_q\end{pmatrix}\ \ \text{and}\ \ J_{n}=\begin{pmatrix} 0& I_n
\\-I_n & 0\end{pmatrix}.$$

\section{Acknowledgements}  The authors are grateful to Martin Svensson
for useful comment on this paper.

\appendix
\section{}\label{app:general}

In this appendix we prove the result stated in Theorem
\ref{theo:rational}.  It shows how the elements of an eigenfamily
$\E$ of complex valued functions on a semi-Riemannian manifold
$(M,g)$ can be used to produce a variety of harmonic morphisms
defined on open and dense subsets of $M$.  The first result shows
how the operators $\tau$ and $\kappa$ behave with respect to
products.

\begin{lemma}\label{lemm:products}
Let $(M,g)$ be a semi-Riemannian manifold and $\E_1,\E_2$ be two
families of complex valued functions on $M$.  If there exist
complex numbers
$\lambda_{1},\mu_{1},\lambda_{2},\mu_{2},\mu\in\cn$ such that for
all $\phi_1,\phi_2\in\E_1$ and $\psi_1,\psi_2\in\E_2$
$$\tau(\phi_1)=\lambda_{1}\phi_1,\ \ \kappa (\phi_1,\phi_2)
=\mu_{1} \phi_1\phi_2,$$ $$\tau(\psi_1)=\lambda_{2}\psi_1,\ \
\kappa (\psi_1,\psi_2)=\mu_{2} \psi_1\psi_2,$$ $$\kappa
(\phi_1,\psi_1)=\mu\phi_1\psi_1$$ then the following relations
hold
$$\tau (\phi_1\psi_1)=(\lambda_{1}+2\mu+\lambda_{2})\phi_1\psi_1,$$
$$\kappa (\phi_1\psi_1,\phi_2\psi_2)
=(\mu_{1}+2\mu+\mu_{2})\phi_1\psi_1\phi_2\psi_2$$
 for all $\phi_1,\phi_2\in\E_1$ and $\psi_1,\psi_2\in\E_2$.
\end{lemma}

\begin{proof}  The statement is an immediate consequence of the
following basic facts concerning first and second order
derivatives of products
$$X(\phi_1\psi_1)=X(\phi_1)\psi_1+\phi_1X(\psi_1),$$
$$X^2(\phi_1\psi_1)=X^2(\phi_1)\psi_1+2X(\phi_1)X(\psi_1)+\phi_1X^2(\psi_1).$$
\end{proof}

The following result shows how the operators $\tau$ and $\kappa$
behave with respect to quotients.

\begin{lemma}\label{lemm:quotient-1}
Let $(M,g)$ be a semi-Riemannian manifold and $P,Q:M\to\cn$ be two
complex valued functions on $M$. If there exists a complex number
$\lambda\in\cn$ such that $$\tau (P)=\lambda P\ \ \text{and}\ \
\tau (Q)=\lambda Q$$ then the quotient $\phi=P/Q$ is a harmonic
morphism if and only if
$$Q^2\kappa (P,P)=PQ\kappa (P,Q)=P^2\kappa (Q,Q).$$
\end{lemma}

\begin{proof}  For first and second order derivatives of the quotient
$P/Q$ we have the following basic facts
$$X(\phi)=\frac{X(P)Q-PX(Q)}{Q^2},$$
$$X^2(\phi)=\frac{Q^2X^2(P)-2QX(P)X(Q)+2PX(Q)X(Q)-PQX^2(Q)}{Q^3}$$
leading to the following formulae for $\tau(\phi)$ and $\kappa
(\phi,\phi )$
$$Q^3\tau(\phi)=Q^2\tau(P)-2Q\kappa(P,Q)+2P\kappa(Q,Q)-PQ\tau(Q),$$
$$Q^4\kappa(\phi,\phi)=Q^2\kappa(P,P)-2PQ\kappa(P,Q)+P^2\kappa(Q,Q).$$
The statement is a direct consequence of these relations.
\end{proof}

\begin{proof}[Proof of Theorem \ref{theo:rational}]
For an eigenfamily $\E=\{\phi_1,\dots ,\phi_n\}$ on the
semi-Riemannian manifold $(M,g)$ we define the infinite sequence
$$\{\E^k\}_{k=1}^\infty$$ by induction $$\E^1=\E\ \ \text{and}
\ \ \E^{k+1}=\E^1\cdot\E^k=\{\phi\cdot\psi |\ \phi\in\E^1,\
\psi\in\E^k\}.$$  It then follows from the fact that $$\kappa
(\phi,\psi )=k\mu_1\phi\psi$$ for all $\phi\in\E^1$ and
$\psi\in\E^k$ and Lemma \ref{lemm:products} that each $\E^{k+1}$
is an eigenfamily on $M$. With this at hand the statement of
Theorem \ref{theo:rational} is an immediate consequence of Lemma
\ref{lemm:quotient-1}.
\end{proof}

\end{document}